\documentclass[12pt]{article}
\usepackage{amsmath}
\usepackage{amsfonts}
\usepackage{amssymb}
\usepackage{graphicx}%
\providecommand{\U}[1]{\protect \rule{.1in}{.1in}}
\newtheorem{theorem}{Theorem}[section]

\newtheorem{definition}[theorem]{Definition}

\newtheorem{lemma}{Lemma}[section]

\newtheorem{proposition}{Proposition}[section]

\newenvironment{proof}[1][Proof]{\noindent \textbf{#1.} }{\  $\Box$}

\title{Reflected BSDE  with stochastic Lipschitz coefficient}

\author{Wen L\"{u}\footnote{Support by the National Basic
Research Program of China (973 Program) grant No. 2007CB814900 and
The Youth Fund of Yantai University (SX08Z9).} \footnote{\emph{Email address:} llcxw@163.com}\\
        \footnotesize{  School of Mathematics, Shandong University, Jinan, \small{250100}, China
 }\\  \footnotesize{School of Mathematics, Yantai University, Yantai 264005, China
 } }

\begin{document}
\date{}
\maketitle
\begin{abstract}
 In this paper,  we deal with a class of
one-dimensional reflected backward stochastic differential equations
with stochastic Lipschitz coefficient.  We derive the existence and
uniqueness of the solutions for those equations via Snell envelope
and the fixed point theorem.
\end{abstract}

\textbf{Keywords:} Reflected backward stochastic differential
equation;  stochastic Lipschitz
 coefficient; Snell envelope

\textbf{AMS 2000  Subject Classification:} 60H10

\section{Introduction}

Nonlinear backward stochastic differential equations (BSDE in short)
were firstly introduced by Pardoux and Peng (1990).  Since then, a
lot of work have been devoted to the study of BSDEs as well as to
their applications. This is due to the connections of BSDEs with
mathematical finance ( see e.g. El Karoui et al. (1997c)) as well as
to stochastic optimal control (see e.g. Peng (1993)) and stochastic
games ( see e.g.  Hamad\`{e}ne and Lepeltier (1995)). El Karoui et
al. (1997a) introduced the notion of one barrier reflected BSDE
(RBSDE in short), which is actually a backward equation but the
solution is forced to stay above a given barrier. This type of BSDEs
is motivated by pricing American options (see El Karoui et al.
(1997b)) and studying the mixed game problems\textbf (see e.g.
Cvitanic and Karatzas (1996), Hamad\`{e}ne and Lepeltier (2000)).

The existence and uniqueness of solution of BSDE in Pardoux and Peng
(1990) and of RBSDE in El Karoui et al. (1997a) are both proved
under the Lipschitz assumption on the coefficient. However, the
Lipschitz condition is too restrictive to be assumed  in many
applications. For instance, the pricing problem of an American claim
is equivalent to solving the linear BSDE
\begin{eqnarray*}
\mbox{d}Y_t=[r(t)Y_t+\theta(t)Z_t]\mbox{d}t+Z_t\mbox{d}B_t,\;
Y_T=\xi,
\end{eqnarray*}
where $r(t)$ is the interest rate and $\theta(t)$ is the risk
premium vector. In general, both of them may be unbounded, therefore
Pardoux and Peng's result may be invalid. And so is it in the case
of RBSDE.

 Consequently, many papers have devoted to
relax the Lipschitz condition in both cases of BSDE and RBSDE (see
e.g. Lepeltier and Martin (1997), El Karoui and Huang (1997), Bender
and Kohlmann (2000), Wang and Huang (2009), Matoussi (1997),
Lepeltier et al. (2005) and the references therein). During them, El
Karoui and Huang (1997) established a general result of existence
and uniqueness for BSDEs driven by a general cadlag martingale  with
stochastic Lipschitz coefficient. Later,  Bender and Kohlmann (2000)
showed the same result for BSDEs driven  by a Brownian motion.
Motivated by the above works, the purpose of the present paper is to
consider a class of one-dimensional RBSDEs with stochastic Lipschitz
coefficient. We try to get the existence and uniqueness of solutions
for those RBSDEs by means of the Snell envelope and the fixed point
theorem.

The rest of the paper is organized as follows.  In Section 2, we
introduce some  notations including some spaces. Section 3 is
devoted to prove the existence and uniqueness of solutions to RBSDEs
with stochastic Lipschitz coefficient.
 \section{Notations}
Let $(B_t)_{t\geq 0}$ be a $d$-dimensional standard Brownian motion
defined on a probability space $(\Omega, \mathcal{F}, {\bf P})$.
 We denote $(\mathcal{F}_t)_{t\geq0}$  the natural filtration of
 $(B_t)_{t\geq 0}$, augmented by all ${\bf P}$-null sets of $\mathcal{F}$.
The Euclidean norm of a vector $y\in {\bf R}^n$ will be
 defined by $|y|$.

Let $T>0$ be a given real number. Let's introduce some spaces:
\begin{itemize}
\item  ${\bf L}^2$ the space of  $\mathcal
{F}_T$-measurable random variables $\xi$ such that
$${\bf E}[|\xi|^2]<+\infty.$$
\item ${\bf S}^2$ the space of predictable processes  $\{\psi_t: t\in[0,
T]\}$ such that $${\bf E}[\sup\limits_{0\leq t\leq
T}|\psi_t|^2]<+\infty.$$
\item ${\bf H}^2$ the space of predictable processes $\{\psi_t: t\in[0, T]\}$ such
that $${\bf E}\int_0^T|\psi_t|^2\mbox{d}t<+\infty.$$
\end{itemize}

Let $\beta>0$ and $(a_t)_{t\geq 0}$ be a nonnegative
$\mathcal{F}_t$-adapted process. Define
 \begin{eqnarray*}
A(t)=\int_0^ta^2(s)\mbox{d}s,\quad 0\leq t\leq T.
\end{eqnarray*}
We further introduce the following spaces:
\begin{itemize}
\item  ${\bf L}^2(\beta, a)$ the space of  $\mathcal
{F}_T$-measurable random variables $\xi$ such that
$${\bf E}[e^{\beta A(T)}|\xi|^2]<+\infty.$$
\item ${\bf S}^2(\beta, a)$ the space of predictable processes  $\{\psi_t: t\in[0,
T]\}$ such that $${\bf E}[e^{\beta A(T)}\sup\limits_{0\leq t\leq
T}|\psi_t|^2]<+\infty.$$
\item ${\bf H}^2(\beta, a)$the space of predictable processes $\{\psi_t: t\in[0, T]\}$ such
that $${\bf E}\int_0^Te^{\beta A(t)}|\psi_t|^2\mbox{d}t<+\infty.$$
\end{itemize}

In this paper, we consider the following RBSDE:
\begin{eqnarray}\label{bsde:1}
\left\{ \begin{array}{l@{\quad \quad}r}Y_t=\xi+\int_t^T f(s, Y_s,
Z_s)\mbox{d}s+K_T-K_t-\int_t^T
Z_s\mbox{d}B_s,\\
Y_t\geq S_t, \; 0\leq t\leq T \;\mbox{a.s. and}\;
\int_0^T(Y_t-S_t)\mbox{d}K_t=0, \mbox{a.s.}
\end{array}\right.
\end{eqnarray}
where the coefficient $f: \Omega\times[0,T]\times {\bf R}\times {\bf R}^d\times {\bf R}\rightarrow {\bf R}$ satisfies the following assumptions:\\
(\textbf{H1}) $\forall t\in[0,T], (y_i,z_i)\in{\bf R}\times{\bf
R}^d, i=1,2$, there are two nonnegative $\mathcal{F}_t$-adapted
processes $\mu(t)$ and $\gamma(t)$ such that
 \begin{eqnarray}\label{lip:1}
|f(t,y_1,z_1)-f(t,y_2,z_2)|\leq \mu(t)|y_1-y_2|+\gamma(t)|z_1-z_2|;
\end{eqnarray} \\
(\textbf{H2})  $\exists\,\epsilon>0$ such that $a^2(t):=
\mu(t)+\gamma^2(t)\geq \epsilon$;\\
(\textbf{H3})  For all $(y, z)\in{\bf R}\times{\bf R}^d$, the
process $f(\cdot,\cdot,y, z)$ is progressively measurable and such
that $\forall t\in[0,T], \frac{\displaystyle f(t,0,0)}{\displaystyle
a}\in {\bf H}^2(\beta, a)$.

 Furthermore, we make the following assumptions:\\
(\textbf{H4}) The terminal value $\xi\in{\bf L}^2(\beta,
a)$;\\
(\textbf{H5}) The "obstacle" $\{S_t, 0\leq t\leq T\}$ is a
continuous progressively measurable real-valued process satisfying
${\bf E}[\sup_{0\leq t\leq T}e^{2\beta A(t)}(S_t^+)^2]<\infty$ and
$S_T\leq\xi$ a.s.

 We now give the
definition of solution to RBSDE (\ref{bsde:1}).
\begin{definition} Let $\beta>0$ and $a$   a nonnegative $\mathcal{F}_t$-adapted
process. A solution to RBSDE (\ref{bsde:1}) is  a triple $(Y,Z,K)$
satisfying (\ref{bsde:1}) such that $(Y,Z)\in {\bf S}^2(\beta,
a)\times {\bf H}^2(\beta, a)$ and $K\in{\bf S}^2$ is continuous and
increasing with $K_0=0$.
\end{definition}

\section{Main results}
\subsection{A priori estimate}
We first give a priori estimate of the solution of RBSDE
(\ref{bsde:1}).
\begin{lemma}\label{lemma:1}
Let $(Y_t,Z_t,K_t)_{0\leq t\leq T}$ be a solution of RBSDE
(\ref{bsde:1}) with data $(\xi, f, T)$. Then there exists a constant
$C_{\beta}$ depending only on $\beta$ such that
\begin{eqnarray*}
&&{\bf E}\left[\sup_{0\leq t\leq T} |Y_t|^2e^{\beta
A(t)}+\int_0^Te^{\beta A(s)}|Z_s|^2\mbox{d}s+\int_0^T e^{\beta
A(s)}a^2(s)|Y_s|^2\mbox{d}s+K_T^2\right]
\\&\leq&
C_{\beta} {\bf E}\left[|\xi|^2e^{\beta A(T)}+\int_0^Te^{\beta
A(s)}\frac{|f(s,0,0)|^2}{ a^2(s)}\mbox{d}s+\sup_{0\leq t\leq
T}e^{2\beta A(t)}(S_t^{+})^2\right].
\end{eqnarray*}
\end{lemma}
\begin{proof} Applying It\^{o}'s formula to $e^{\beta A(t)}|Y_t|^2$, we have
\begin{eqnarray*}
&&e^{\beta A(t)}|Y_t|^2+\int_t^Te^{\beta A(s)}|Z_s|^2\mbox{d}s+\beta
\int_t^Ta^2(s)e^{\beta A(s)}|Y_s|^2\mbox{d}s
\\&=&
e^{\beta A(T)}|\xi|^2+2\int_t^Te^{\beta
A(s)}Y_sf(s,Y_s,Z_s)\mbox{d}s+2\int_t^Te^{\beta
A(s)}Y_s\mbox{d}K_s-2\int_t^Te^{\beta A(s)}Y_sZ_s\mbox{d}B_s
\\&\leq&
e^{\beta A(T)}|\xi|^2+\frac{\beta}{2}\int_t^T a^2(s)e^{
A(s)}|Y_s|^2\mbox{d}s+2\int_t^Te^{\beta A(s)}\frac{1}{\beta
a^2(s)}|f(s,Y_s,Z_s)|^2\mbox{d}s
\\&&
+2\int_t^Te^{\beta A(s)}Y_s\mbox{d}K_s-2\int_t^Te^{\beta
A(s)}Y_sZ_s\mbox{d}B_s
\\&\leq& e^{\beta
A(T)}|\xi|^2+\frac{\beta}{2}\int_t^T a^2(s)e^{
A(s)}|Y_s|^2\mbox{d}s+\frac{6}{\beta}[\int_t^Te^{\beta
A(s)}a^2(s)|Y_s|^2\mbox{d}s+\int_t^Te^{\beta A(s)}|Z_s|^2\mbox{d}s]
\\&& +\frac{6}{\beta}\int_t^Te^{\beta
A(s)}\frac{|f(s,0,0)|^2}{a^2(s)}\mbox{d}s +2\int_t^Te^{\beta
A(s)}Y_s\mbox{d}K_s-2\int_t^Te^{\beta A(s)}Y_sZ_s\mbox{d}B_s.
\end{eqnarray*}
Consequently,
\begin{eqnarray}
e^{\beta A(t)}|Y_t|^2&+&(1-\frac{6}{\beta})\int_t^Te^{\beta
A(s)}|Z_s|^2\mbox{d}s+(\frac{\beta}{2}-\frac{6}{\beta})
\int_t^Ta^2(s)e^{\beta A(s)}|Y_s|^2\mbox{d}s\nonumber
 \nonumber\\&\leq&
  e^{\beta
A(T)}|\xi|^2 +\frac{6}{\beta}\int_t^Te^{\beta
A(s)}\frac{|f(s,0,0)|^2}{a^2(s)}\mbox{d}s \nonumber
\\&&\label{estam:2}
+2\int_t^Te^{\beta A(s)}S_s\mbox{d}K_s-2\int_t^Te^{\beta
A(s)}Y_sZ_s\mbox{d}B_s.
\end{eqnarray}
where we have used the fact that $\mbox{d}K_s={\bf
I}_{[Y_s=S_s]}\mbox{d}K_s$  and the stochastic Lipschitz property of
$f$. For a sufficient large $\beta>0$, taking
 expectation on both sides above, we get
\begin{eqnarray}\label{estam:3}
&&{\bf E}[\int_t^T a^2(s)e^{\beta
A(s)}|Y_s|^2\mbox{d}s+\int_t^Te^{\beta
A(s)}|Z_s|^2\mbox{d}s]\nonumber
\\&\leq&
c_{\beta}{\bf E}\left[e^{\beta A(T)}|\xi|^2+\int_t^Te^{\beta
A(s)}\frac{|f(s,0,0)|^2}{a^2(s)}\mbox{d}s +2\int_t^Te^{\beta
A(s)}S_s^+\mbox{d}K_s\right].
\end{eqnarray}
Moreover, by the Burkholder-Davis-Gundy's inequality, one can derive
that
\begin{eqnarray*}
&&{\bf E}[\sup_{0\leq t\leq T}|\int_t^Te^{\beta
A(s)}Y_sZ_s\mbox{d}B_s|]
\\&\leq&
{\bf E}[|\int_0^Te^{\beta A(s)}Y_sZ_s\mbox{d}B_s|]+{\bf
E}[\sup_{0\leq t\leq T}|\int_0^te^{\beta A(s)}Y_sZ_s\mbox{d}B_s|]
\\&\leq&
2c{\bf E}\left\{\left[\int_0^Te^{2\beta
A(s)}|Y_s|^2|Z_s|^2\mbox{d}s\right]^{\frac{1}{2}}\right\}
\\&\leq&
2c{\bf E} \left[(\sup_{0\leq t\leq T}e^{\beta
A(t)}|Y_t|^2)^{\frac{1}{2}}(\int_0^Te^{\beta
A(s)}|Z_s|^2\mbox{d}s)^{\frac{1}{2}}\right]
\\&\leq&
\frac{1}{2}{\bf E}[(\sup_{0\leq t\leq T}e^{\beta
A(t)}|Y_t|^2)]+2c^2{\bf E}[\int_0^Te^{\beta A(s)}|Z_s|^2\mbox{d}s].
\end{eqnarray*}
Combining this with (\ref{estam:2}) and (\ref{estam:3}), we have
\begin{eqnarray}\label{estam:4}
&&{\bf E}\left[\sup_{0\leq t\leq T} |Y_t|^2e^{\beta A(t)}+\int_0^T
a^2(s)e^{\beta A(s)}|Y_s|^2\mbox{d}s+\int_0^Te^{\beta
A(s)}|Z_s|^2\mbox{d}s\right]\nonumber
\\&\leq&k_{\beta}{\bf E}\left[e^{\beta
A(T)}|\xi|^2+\int_0^Te^{\beta
A(s)}\frac{|f(s,0,0)|^2}{a^2(s)}\mbox{d}s +2\int_0^Te^{\beta
A(s)}S_s^+\mbox{d}K_s\right].
\end{eqnarray}

We now give an  estimate of $K_T^2$. From the equation
\begin{eqnarray*}
K_T=Y_0-\xi-\int_0^Tf(s,Y_s,Z_s)\mbox{d}s+\int_0^TZ_s\mbox{d}B_s
\end{eqnarray*}
and (\ref{estam:4}), we have
\begin{eqnarray*}
&&{\bf E}[K_T^2] \\&\leq&
 d_{\beta}{\bf E}\left[\sup_{0\leq t\leq T}
|Y_t|^2e^{\beta A(t)}+|\xi|^2 +\int_0^T|Z_s|^2\mbox{d}s\right.
\\&&\left.+\int_0^Ta^2(s)e^{-\beta A(s)}\mbox{d}s\int_0^Te^{\beta
A(s)}\frac{|f(s,Y_s,Z_s)|^2}{a^2(s)}\mbox{d}s\right]
\\
&\leq&
 d_{\beta}{\bf E}\left[e^{\beta A(T)}|\xi|^2+\int_0^Ta^2(s)e^{\beta
A(s)}|Y_s|^2\mbox{d}s +\int_0^Te^{\beta A(s)}|Z_s|^2\mbox{d}s\right.
\\&&\left.+\int_0^Te^{\beta
A(s)}\frac{|f(s,0,0)|^2}{a^2(s)}\mbox{d}s\right]
\\
&\leq& d_{\beta}{\bf E}\left[e^{\beta A(T)}|\xi|^2+\int_0^Te^{\beta
A(s)}\frac{|f(s,0,0)|^2}{a^2(s)}\mbox{d}s +2\int_0^Te^{\beta
A(s)}S_s^+\mbox{d}K_s\right]
\\
&\leq& d_{\beta}{\bf E}\left[e^{\beta A(T)}|\xi|^2+\int_0^Te^{\beta
A(s)}\frac{|f(s,0,0)|^2}{a^2(s)}\mbox{d}s+\sup_{0\leq t\leq T}
e^{2\beta A(t)}(S_t^+)^2\right]+\frac{1}{2}{\bf E}[K_T^2].
\end{eqnarray*}
Hence,
\begin{eqnarray}\label{estam:5}
{\bf E}[K_T^2] \leq d_{\beta}{\bf E}\left[e^{\beta
A(T)}|\xi|^2+\int_0^Te^{\beta
A(s)}\frac{|f(s,0,0)|^2}{a^2(s)}\mbox{d}s+\sup_{0\leq t\leq T}
e^{2\beta A(t)}(S_t^+)^2\right],
\end{eqnarray}
where we use the notation $d_{\beta}$ for a constant depending only
on $\beta$ and whose value could be changing from line to line. We
get the desired  result by estimates (\ref{estam:5}) and
(\ref{estam:4}).
\end{proof}
\subsection{Existence and uniqueness of solution}
We first consider the special case that is the  coefficient does not
depend on $(Y,Z)$, i.e. $f(\omega,t, y,z)\equiv g(\omega, t)$. We
have the following result.

\begin{theorem}\label{theorem:1}
Let  $\beta>0$ large enough and $a$ a nonnegative
$\mathcal{F}_t$-adapted process. Assume $\frac{\displaystyle
g}{\displaystyle a}\in{\bf H}^2(\beta, a)$ and
(\textbf{H4})-(\textbf{H5}) hold. Then RBSDE (\ref{bsde:1}) with
data $(\xi, g, S)$  has a solution.
\end{theorem}
\begin{proof}
 For $0\leq t\leq T$, we define
\begin{eqnarray*}
\widetilde{Y}_t=\mbox{ess}\sup\limits_{\nu\geq t}{\bf
E}[\int_0^{\nu}g(s)ds +S_{\nu}{\bf I}_{\{\nu<T\}}+\xi{\bf
I}_{\{\nu=T\}}|\mathcal {F}_{t}]
\end{eqnarray*}
where $\nu$ is a $\mathcal{F}_t$-stopping time. The process
$\widetilde{Y}_t$ is called the Snell envelope of the process which
is inside $\mbox{ess}\sup$.

By assumptions of the theorem, one can easily  to check that
$\xi\in{\bf L}^2$, $S_t^+\in{\bf S}^2$  and
$(\int_0^t|g(s)|\mbox{d}s)_{0\leq t\leq T}\in{\bf L}^2$. Indeed, for
given $\beta>0$, by H\"{o}lder inequality, we have
\begin{eqnarray*}
{\bf E}\left[\left(\int_0^T|g(s)|\mbox{d}s\right)^2\right]&=&{\bf
E}\left[\left(\int_0^T|\frac{g(s)}{a(s)}||a(s)|\mbox{d}s\right)^2\right]
\\&\leq&{\bf E}\left[\left(\int_0^T |\frac{g(s)}{a(s)}|^2e^{\beta A(s)}\mbox{d}s\right)\left(\int_0^T a(s)^2e^{-\beta
A(s)}\mbox{d}s\right)\right]
\\&\leq&\frac{1}{\beta}{\bf E}\left[\left(\int_0^T|\frac{g(s)}{a(s)}|^2e^{\beta A(s)}\mbox{d}s\right)\right]<+\infty.
\end{eqnarray*}
Consequently, by Doob-Meyer decomposition theorem in Dellacherie and
Meyer (1980), there exists an increasing continuous process
$(K_t)_{0\leq t\leq T}$ which belongs to ${\bf S}^2$ ($K_0=0$) and a
martingale $M_t\in{\bf S}^2$ such that
\begin{eqnarray*}
 \widetilde{Y}_t=M_t-K_t.\quad\forall   t\in[0, T]
\end{eqnarray*}
Since $M_t\in {\bf S}^2$, there exists $Z_t\in {\bf H}^2$ such that
\begin{eqnarray*}
 M_t=M_0+\int_0^tZ_s\mbox{d}B_s.\quad\forall t\in[0, T]
\end{eqnarray*}
Let $Y_t=\widetilde{Y}_t-\int_0^tf(s)\mbox{d}s$, by Proposition 5.1
of El Karoui et al. (1997a), we derive that the triple $(Y,Z,K)$
verities
\begin{eqnarray*}
Y_t=\xi+\int_t^Tg(s)\mbox{d}s+K_T-K_t-\int_t^TZ_s\mbox{d}B_s.
\end{eqnarray*}
Moreover, $Y_t\geq S_t$ and $\int_0^T(Y_t-S_t)\mbox{d}K_t=0$. By
Lemma \ref{lemma:1}, $(Y_t,Z_t,K_t)_{0\leq t\leq T}$ is a solution
of RBSDE (\ref{bsde:1}).
\end{proof}

Furthermore, we have the following uniqueness result.
\begin{proposition}\label{proposition:1}
With the same assumptions of Theorem \ref{theorem:1}, the RBSDE
(\ref{bsde:1}) with data  $(\xi, g, S)$  has at most one solution.
\end{proposition}
{\bf Proof.} Let $(Y,Z,K)$ and $(Y',Z',K')$ be two solutions of
RBSDE (\ref{bsde:1}). Let
\begin{eqnarray*}
\Delta Y=Y-Y', \;\Delta Z=Z-Z',\;\Delta K=K-K'.
\end{eqnarray*}
For $0\leq t\leq T$, we  have
\begin{eqnarray*}
\Delta Y_t=\Delta K_T-\Delta K_t-\int_t^T\Delta Z_s\mbox{d}B_s.
\end{eqnarray*}
 Applying It\^{o}'s formula to $e^{\beta A(t)}|\Delta
Y_t|^2$, we obtain
\begin{eqnarray*}
-{\bf E}[e^{\beta A(t)}|\Delta Y_t|^2]=-2\,{\bf E}[\int_t^Te^{\beta
A(s)}\Delta Y_s\mbox{d}(\Delta K_s)]+{\bf E}[\int_t^Te^{\beta
A(s)}|\Delta Z_s|^2\mbox{d}s].
\end{eqnarray*}
Noting that $\int_t^Te^{\beta A(s)}\Delta Y_s\mbox{d}(\Delta
K_s)\leq 0$, it follows that $\Delta Y_t=\Delta Z_t=0$ and then
$\Delta K_t=0$, $0\leq t\leq T$ a.s.\quad $\Box$

\bigskip
We can now  state and prove our main result.
\begin{theorem}\label{theorem:4}
Assume (H1)-(H5) hold for a sufficient large $\beta$. Then  RBSDE
(\ref{bsde:1}) with data $(\xi, f, S)$ has a unique solution.
\end{theorem}
\begin{proof}
 Let $\mathcal{H}(\beta, a)={\bf S}^2(\beta, a)\times{\bf H}^2(\beta, a)$. Given $(U,V)\in\mathcal{H}(\beta, a)$,
consider the following RBSDE:
\begin{eqnarray}\label{bsde:2}
Y_t=\xi+\int_t^Tf(s, U_s,
V_s)\mbox{d}s+K_T-K_t-\int_t^TZ_s\mbox{d}B_s.
\end{eqnarray}
By Young's inequality, we have
\begin{eqnarray*}
\frac{|f(t,U_t,V_t)|^2}{a^2(t)}\leq
3[a^2(t)|U_t|^2+|V_t|^2+\frac{|f(t,0,0)|^2}{a^2(t)}],
\end{eqnarray*}
it follows from (H3) and Theorem \ref{theorem:1} that the RBSDE
(\ref{bsde:2}) has a unique solution.

Define a mapping $\Phi$ from $\mathcal{H}(\beta, a)$ to itself. Let
$(U',V')$ be another element in $\mathcal{H}(\beta, a)$, set
\begin{eqnarray*}
(Y,Z)=\Phi(U,V),\;(Y',Z')=\Phi(U',V'),
\end{eqnarray*}
where $(Y,Z,K)$ (resp. $(Y', Z',K')$) is the unique solution of the
RBSDE associated with data $(\xi, f(t, U_t, V_t), S)$ (resp.$(\xi,
f(t, U'_t, V'_t),S)$).

Let
\begin{eqnarray*}
&\Delta Y=Y-Y', \Delta Z=Z-Z', \Delta U=U-U', \Delta V=V-V',
\\&\Delta f_s=f(s, U_s, V_s)-f(s, U'_s, V'_s), \Delta K=K-K'.
\end{eqnarray*}
For $0\leq t\leq T$, we  have
\begin{eqnarray*}
\Delta Y_t= \int_t^T\Delta f_s\mbox{d}s+ \Delta K_T- \Delta
K_t-\int_t^T \Delta Z_s\mbox{d}B_s
\end{eqnarray*}
Applying It\^{o}'s formula to $e^{\beta A(t)}|\Delta Y_t|^2$, using
(\textbf{H1}) and the fact that $\mbox{d}K_s={\bf
I}_{[Y_s=S_s]}\mbox{d}K_s$ and $\mbox{d}K'_s={\bf
I}_{[Y'_s=S_s]}\mbox{d}K'_s$, we get
\begin{eqnarray*}
&&e^{\beta A(t)}|\Delta Y_t|^2+\beta\int_t^T a(s)^2e^{\beta
A(s)}|\Delta Y_s|^2\mbox{d}s+\int_t^T e^{\beta A(s)}|\Delta
Z_s|^2\mbox{d}s
\\&\leq&
2\int_t^Te^{\beta A(s)}\Delta Y_s\Delta f_s\mbox{d}s
+2\int_t^Te^{\beta A(s)} \Delta Y_s \mbox{d}(\Delta K_s)-\int_t^T
2e^{\beta A(s)} \Delta Y_s \Delta Z_s\mbox{d}B_s
\\&\leq&
 2\int_t^Te^{\beta A(s)} \Delta Y_s \Delta f_s \mbox{d}s
 -\int_t^T 2e^{\beta A(s)} \Delta Y_s  \Delta Z_s\mbox{d}B_s
\\&\leq& \frac{\beta}{2}\int_t^T a(s)^2e^{\beta A(s)}|\Delta Y_s|^2\mbox{d}s
+\frac{6}{\beta}\int_t^Te^{\beta A(s)}| ( a(s)^2|\Delta
U_s|^2+|\Delta V|^2)\mbox{d}s\\&&-\int_t^T 2e^{\beta A(s)} \Delta
Y_s \Delta Z_s\mbox{d}B_s,
\end{eqnarray*}
it follows that
\begin{eqnarray*}
&&{\bf E}[\int_t^T a(s)^2e^{\beta A(s)}|\Delta Y_s|^2\mbox{d}s]+{\bf
E}[\int_t^T e^{\beta A(s)}|\Delta
Z_s|^2\mbox{d}s]\\&\leq&(\frac{12}{\beta^2}+\frac{6}{\beta})\left\{{\bf
E}[\int_t^T a(s)^2e^{\beta A(s)}|\Delta U_s|^2]+{\bf
E}[\int_t^Te^{\beta A(s)}|\Delta V|^2\mbox{d}s]\right\}.
\end{eqnarray*}
For $\beta>0$ large enough, one can easily to check that $\Phi$ is a
contraction mapping with the norm
\begin{eqnarray*}
\|(Y,Z)\|^2_{\beta}={\bf E}\left[\int_0^T e^{\beta
A(s)}(a(s)^2|Y_s|^2+|Z_s|^2)ds\right].
\end{eqnarray*}
Thus, $\phi$ has a unique fixed point and the theorem is
proved.
\end{proof}

\end{document}